\documentclass{article}
\usepackage{amssymb}
\usepackage{amsmath}
\usepackage{abstract}
\usepackage{indentfirst}
\usepackage[utf8]{inputenc}
\usepackage[hidelinks]{hyperref}

\usepackage{titlesec}
\titleformat{\section}[block]{\center\scshape\large}{\thesection.}{0.7em}{}
\titleformat{\subsection}[block]{\center}{\thesubsection.}{0.7em}{}
\renewcommand\thesection{\Roman{section}}

\usepackage{amsthm}
\newtheorem{theorem}{Theorem}
\newtheorem{lemma}{Lemma}

\newtheorem{proposition}{Proposition}
\theoremstyle{definition}

\newcommand{\keywords}[1]
{
	\small
	\textbf{Keywords:} #1
}

\title{A Corrected Simplified Proof of Chen's Theorem}

\usepackage{authblk}
\author[a]{Zihao Liu}
\affil[a]{International Department, The Affiliated High School of SCNU,\authorcr Email: \url{mailto:travor_lzh@163.com}}

\date{}

\begin{document}
	\maketitle
	\begin{abstract}
		In 1973, J.-R. Chen \cite{chen_representation_1973} showed that every large even integer is a sum of a prime and a product of at most two primes. In this paper, the author indicates and fixes the issues in a simplified proof of this result given by Pan et al. \cite{pan_representation_1975}
	\end{abstract}
	\section{Introduction}
	For brevity, let $\{1,c\}$ denote the following proposition:
	
	There exists a constant $N_0$ such that for all even $N>N_0$, there exists a prime $p<N-1$ such that $N-p$ is a product of at most $c$ primes.
	
	Initially, the study of $\{1,c\}$ assumes the truthfulness of the Generalized Riemann Hypothesis for Dirichlet L-function:

	\begin{proposition}[GRH]
		Let $\pi(x;q,l)$ denote the number of primes $p$ such that
		\[
			p\le x,p\equiv l\pmod q
		\]
		If $(l,q)=1$, then
		\begin{equation}
			\label{eqn-grh}
			\pi(x;q,l)={\operatorname{li}x\over\varphi(q)}+O(\sqrt x\log x)
		\end{equation}
		where $\operatorname{li}x$ denotes the logarithmic integral function.
	\end{proposition}
	
	Assuming GRH, Estermann \cite{estermann_neue_1932} proves $\{1,6\}$ in 1932. Subsequently, this result is improved to $\{1,4\}$ by Wang \cite{wang_representation_1956} using the sieves of Buchstab \cite{wang_goldbach_2002} and Selberg \cite{selberg_elementary_1947}. Then in 1962, Wang \cite{wang_representation_1962} applies Kuhn's methods \cite{kuhn_uber_1954} and shows $\{1,3\}$ holds under GRH.

	In 1965, Bombieri \cite{bombieri_large_1965} and Vinogradov \cite{vinogradov_uber_1965} apply the large sieve method \cite{davenport_multiplicative_1980} and prove the following unconditional mean value theorem:

	\begin{theorem}[Bombieri-Vinogradov]
		\label{th-bv}
		For all fixed $A>0$, there exists $B=B(A)$ such that for all large $x$:
		\[
			\sum_{q\le x^{\frac12}\log^{-B}x}\max_{y\le x}\max_{(l,q)=1}\left|\pi(x;q,l)-{\operatorname{li}x\over\varphi(q)}\right|\ll{x\over\log^Ax}
		\]
	\end{theorem}

	By replacing \eqref{eqn-grh} with \autoref{th-bv} while estimating error terms, Estermann's $\{1,6\}$ and Wang's $\{1,4\}$ and $\{1,3\}$ become unconditional results.

	To prove $\{1,2\}$, Chen applies the Jurkat--Richert theorem\footnote{Chen himself regards this as Richert's sieve \cite{richert_selbergs_1969}, but the sieve method in Richert's paper is essentially the Jurkat--Richert theorem \cite{jurkat_improvement_1965}.} and proves\footnote{See Lemma 9 of \cite{chen_representation_1973}} that if $W(N)$ denotes the number of primes $p$ such that $N-p$ has no prime factors $\le N^{1\over10}$ and has most one prime factor in $(N^{1\over10},N^{\frac13}]$, then for sufficiently large even $N$

	\begin{equation}
		\label{eqn-wn}
		W(N)\ge2.6408\mathfrak S(N){N\over\log^2N}
	\end{equation}

	where $\mathfrak S(N)$, the singular series for the Goldbach's problem, is defined as

	\begin{equation*}
		\mathfrak S(N)=\prod_{\substack{p|N\\p>2}}{p-1\over p-2}\prod_{p>2}\left(1-{1\over(p-1)^2}\right)
	\end{equation*}

	It follows from pigeonhole principle that $W(N)$ gives a lower bound for the number of ways to express $N$ as a sum of a prime and a product of at most three primes, and $N-p$ counted by $W(N)$ has exactly three prime factors if and only if

	\begin{equation}
		\label{eqn-p3}
		N-p=p_1p_2p_3,N^{1\over10}<p_1\le N^{\frac13}<p_2<p_3
	\end{equation}

	To estimate primes satsifying \eqref{eqn-p3}, Chen devised the following switched sum

	\begin{equation*}
		\Omega=\sum_a\sum_{\substack{ap_3\le N\\N-ap_3\text{ prime}}}f(a)
	\end{equation*}
	where $f(a)$ is the characteristic function for the condition that
	\begin{equation}
		\label{eqn-adef}
		a=p_1p_2,N^{1\over10}<p_1\le N^{\frac13}<p_2\le(N/p_1)^{\frac12} \\
	\end{equation}

	Then, he applies the Selberg's sieve \cite{selberg_elementary_1947} and the large sieve to obtain

	\begin{equation}
		\label{eqn-omegabound}
		\Omega\le3.9404\mathfrak S(N){N\over\log^2N}
	\end{equation}

	By symmetry, it is evident that the number of primes $p$ satisfying \eqref{eqn-p3} is bounded by half of $\Omega$, so $\{1,2\}$ is deduced as the number of ways to express $N$ as a sum of a prime and a product of at most two primes is no less than

	\[
		W(N)-\frac\Omega2\ge\left(2.6408-{3.9404\over2}\right)\mathfrak S(N){N\over\log^2N}>0.67\mathfrak S(N){N\over\log^2N}
	\]

	Soon after Chen, Pan et al. \cite{pan_representation_1975} simplified the proof of \eqref{eqn-omegabound} by introducing the following new mean value theorem:

	\begin{theorem}[Pan et al.]
		\label{th-mvt}
		Let $\pi(x;a,q,l)$ denote the number of primes satisfying
		\[
			ap\le x,ap\equiv l\pmod q
		\]
		and $\Delta(x;a,q,l)$ be defined by
		\begin{equation}
			\label{eqn-delta-def}
			\Delta(x;a,q,l)=\pi(x;a,q,l)-{\operatorname{li}\frac xa\over\varphi(q)}
		\end{equation}
		Then for every fixed $A>0$, there exists $B=B(A)$ such that
		\begin{equation*}
			\sum_{q\le x^{\frac12}\log^{-B}x}\mu^2(q)3^{\omega(q)}\max_{y\le x}\max_{(l,q)=1}\left|\sum_{(a,q)=1}f(a)\Delta(y;a,q,l)\right|\ll{x\over\log^Ax}
		\end{equation*}
	\end{theorem}

	\autoref{th-mvt} is devoted to give upper estimate for the error terms emerged from the estimation of $\Omega$. However, this result does not serve to estimate all the error terms appeared in the upper bound sieve for $\Omega$.

	In this paper, we present a corrected version of Pan et al.'s proof of \eqref{eqn-omegabound}. The issues in their original proof are discussed and resolved in \autoref{sn-r}.
	
	\section{Auxiliary lemmas}
	\begin{lemma}
		\label{lm-mb}
		Let $\beta>\alpha>0$ be fixed. Then for all $x\ge2$, the sum
		\[
			\sum_{x^\alpha<p\le x^\beta}\frac1p
		\]
		is bounded.
	\end{lemma}
	\begin{proof}
		Mertens' second theorem asserts that there exists a fixed constant $B_1$ such that
		\begin{equation}
			\label{eqn-mertens}
			\sum_{p\le x}\frac1p=\log\log x+B_1+O\left(1\over\log x\right)
		\end{equation}
		so we have
		\[
			\sum_{x^\alpha<p\le x^\beta}\frac1p=\log\frac\beta\alpha+O\left(1\over\log x\right)
		\]
		This suggests the sum is bounded.
	\end{proof}
	\begin{lemma}
		\label{lm-muomega}
		Let $A>0$ be a fixed constant. Then for any integer $n\ge1$
		\[
			\sum_{d|n}{\mu^2(d)A^{\omega(d)}\over\varphi(d)}\ll(\log\log3n)^A
		\]
		where $\omega(n)$ denotes the number of distinct prime factors of $n$. In particular, if $n$ is the product of distinct primes $\le y$, then
		\begin{equation}
			\label{eqn-muomega}
			\sum_{d|n}{\mu^2(d)A^{\omega(d)}\over\varphi(d)}=\sum_{d|n}{A^{\omega(d)}\over\varphi(d)}\ll(\log y)^A
		\end{equation}
	\end{lemma}
	\begin{proof}
		It is evident that the left hand side is multiplicative, so by definition we have
		\[
			\sum_{d|n}{\mu^2(d)A^{\omega(d)}\over\varphi(d)}=\prod_{p|n}\left(1+{A\over p-2}\right)\le\exp\left\{\sum_{p|n}{A\over p-1}\right\}
		\]
		To estimate the remaining sum, we replace the denominator:
		\[
			\sum_{p|n}{1\over p-1}=\sum_{p|n}\frac1p+\sum_{p|n}{1\over p(p-1)}=\sum_{p|n}\frac1p+O(1)
		\]
		Now, we introduce a parameter $1\le u\le n$ so that
		\begin{align*}
			\sum_{p|n}\frac1p
			&\le\sum_{p\le u}\frac1p+\sum_{\substack{p|n\\p>u}}\frac1p
			\le\sum_{p\le u}\frac1p+{\omega(n)\over u} \\
			&=\log\log u+O(1)+O\left(\log n\over u\right)
		\end{align*}
		where the last line follows from \eqref{eqn-mertens} and the fact that $\omega(n)=O(\log n)$. Plugging this result back with $u=\log 3n$, we obtain
		\[
			\sum_{d|n}{\mu^2(d)A^{\omega(d)}\over\varphi(d)}\le\exp\{A\log\log\log 3n+O(1)\}\ll(\log\log 3n)^A
		\]
		If $n$ is the product of primes $\le y$, then by Chebyshev's estimates we have
		\[
			\log n=\sum_{p\le y}\log p\ll y
		\]
		Thus, \eqref{eqn-muomega} is also proven.
	\end{proof}
	\begin{lemma}[Selberg's sieve]
		\label{lm-selberg}
		Let $Q$ denote the product of primes $\le z=N^{\frac14-\frac\varepsilon2}$ that do not divide $N$.
		There exists real sequence $\lambda_d$ satisfying
		\begin{enumerate}
			\item $\lambda_1=1$.
			\item $\lambda_d=0$ for $d>z$ or $d\nmid Q$.
			\item $|\lambda_d|\le1$.
		\end{enumerate}
		such that
		\begin{equation*}
			\sum_{d_1,d_2}{\lambda_{d_1}\lambda_{d_2}\over\varphi([d_1,d_2])}=[8+O(\varepsilon)]{\mathfrak S(N)\over\log N}
		\end{equation*}
	\end{lemma}
	\begin{proof}
		See Lemma 2 and §4 of \cite{wang_representation_1962}.
	\end{proof}
	\begin{lemma}
		\label{lm-mt}
		For large $N$, we have
		\begin{equation*}
			\sum_a{f(a)\over a\log(N/a)}\le{0.49254\over\log N}
		\end{equation*}
	\end{lemma}
	\begin{proof}
		According \autoref{eqn-adef} and \autoref{eqn-mertens}, we can use partial summation to transform the leftmost sum into

		\begin{align*}
			&=\sum_{N^{1\over10}<p_1\le N^{\frac13}}\sum_{N^{\frac13}<p_2\le(N/p_1)^{\frac12}}{1\over p_1p_2\log(x/p_1p_2)} \\
			&\sim\left.\int_{N^{1\over10}}^{N^{\frac13}}{\mathrm dx\over x\log x}\int_{N^{\frac13}}^{(N/x)^{\frac12}}{1\over\log N-\log x-\log y}{\mathrm dy\over y\log y}\right\rbrace\substack{\alpha=\log x/\log N\\\beta=\log y/\log N} \\
			&={1\over\log N}\int_{1\over10}^{\frac13}{\mathrm d\alpha\over\alpha}\int_{\frac13}^{{1-\alpha\over2}}{\mathrm d\beta\over\beta(1-\alpha-\beta)}<{0.49254\over\log N}
		\end{align*}
		where the last inequality follows the numerical calculations in (28) of \cite{chen_representation_1973}.
	\end{proof}
	\section{Preliminary treatments for $\Omega$}
	Since every prime number is either co-prime to $Q$ or a divisor of $Q$, we have
	\begin{align*}
		\Omega
		&\le\sum_a\sum_{\substack{ap\le N\\(N-ap,Q)=1}}f(a)+N^{\frac23}z \\
		&\le\underbrace{\sum_af(a)\sum_{ap\le N}\left(\sum_{d|(N-ap)}\lambda_d\right)^2}_M+N^{11\over12}
	\end{align*}
	where $\lambda_d$ is defined as in \autoref{lm-selberg}. By an interchanging of summation, there is
	\begin{align*}
		M
		&=\sum_af(a)\sum_{ap\le N}\sum_{d_1,d_2|(N-ap,Q)}\lambda_{d_1}\lambda_{d_2} \\
		&=\sum_{d_1,d_2}\lambda_{d_1}\lambda_{d_2}\sum_a\sum_{\substack{ap\le N\\ [d_1,d_2]|(N-ap)}}f(a) \\
		&=\sum_{d_1,d_2}\lambda_{d_1}\lambda_{d_2}\sum_af(a)\pi(N;a,[d_1,d_2],N)=M_1+R
	\end{align*}
	where $M_1$ and $R$ satisfy:
	\begin{equation}
		\label{eqn-m1}
		M_1=\sum_{d_1,d_2}{\lambda_{d_1}\lambda_{d_2}\over\varphi([d_1,d_2])}\sum_af(a)\operatorname{li}\frac Na
	\end{equation}
	\begin{equation}
		\label{eqn-r}
		|R|\le\sum_{\substack{d|Q\\ d\le N^{\frac12-\varepsilon}}}3^{\omega(d)}\left|\sum_af(a)\Delta(N;a,d,N)\right|
	\end{equation}
	where the factor $3^{\omega(d)}$ follows from the fact that there are exactly $3^{\omega(d)}$ pairs of $d_1,d_2$ satisfying $[d_1,d_2]=d$. To estimate the main term, notice that
	\[
		\operatorname{li}x={x\over\log x}+O\left(x\over\log^2x\right)
	\]
	so plugging \autoref{lm-selberg} and \autoref{lm-mt} into \eqref{eqn-m1} gives

	\begin{equation}
		\label{eqn-m10}
		M_1\le[8+O(\varepsilon)]{\mathfrak S(N)\over\log N}\cdot{0.49254N\over\log N}\le3.94033\mathfrak S(N){N\over\log^2N}
	\end{equation}

	Thus the remaining task is to estimate the error term in \eqref{eqn-r}.

	\section{Salvaging the estimation of $R$}\label{sn-r}
	
	The authors of \cite{pan_representation_1975} concluded
	\begin{equation}
		\label{eqn-rbound}
		R\ll{N\over\log^AN}
	\end{equation}
	merely from \autoref{th-mvt} because they assume $f(a)\ne0$ automatically ensures $a$ is co-prime to $d$. However, due to \eqref{eqn-adef}, it is possible for $a$ to possess a prime divisor $\le N^{\frac14-\frac\varepsilon2}$, so their arguments for \eqref{eqn-rbound} are incorrect.

	To salvage \eqref{eqn-rbound} and their proof of Chen's theorem, we consider each case separately. After applying \autoref{th-mvt} to estimate the sum over $(a,d)=1$, all we need is to estimate the remaining situation where $(a,d)>1$:
	\begin{equation}
		\label{eqn-r0}
		R\ll{N\over\log^AN}+\underbrace{\sum_{\substack{d|Q\\d\le N^{1/2}}}3^{\omega(d)}\sum_{(a,d)>1}f(a)|\Delta(N;a,d,N)|}_{R_1}
	\end{equation}
	Since $(a,d)>1$ implies $\pi(N;a,d,N)\le1$, it follows from \eqref{eqn-delta-def} and \autoref{lm-muomega} that
	\begin{align*}
		R_1
		&\ll\sum_{d|Q}{3^\omega(d)\over\varphi(d)}\max_{d\le N^{1/2}}\sum_{(a,d)>1}f(a)\operatorname{li}\frac Na \\
		&\ll(\log N)^3\max_{d\le N^{1/2}}\sum_{(a,d)>1}f(a){N\over a\log(N/a)} \\
		&\ll N(\log N)^2\max_{d\le N^{1/2}}\sum_{(a,d)>1}{f(a)\over a}
	\end{align*}
	where the last $\ll$ follows from the fact that $f(a)\ne0$ implies $a\le N^{\frac23}$. Due to \eqref{eqn-adef}, $(a,d)>1$ if and only if $p_1|d$, so
	\begin{align*}
		\sum_{(a,d)>1}{f(a)\over a}
		&=\sum_{\substack{N^{1\over10}<p_1\le N^{\frac13}\\ p_1|d}}{1\over p_1}\sum_{N^{\frac13}<p_2\le(N/p_1)^{\frac12}}{1\over p_2} \\
		&\le\sum_{\substack{p_1>N^{1\over10}\\ p_1|d}}{1\over p_1}\sum_{N^{\frac13}<p_2\le N^{4\over10}}{1\over p_2}
	\end{align*}

	Now, we apply \autoref{lm-mb} with $x=N,\alpha=\frac13,\beta={4\over10}$, so that
	\[
		\sum_{(a,d)>1}{f(a)\over a}\ll\sum_{\substack{p>N^{1\over10}\\p|d}}\frac1p<{10N^{-{1\over10}}\over\log N}\sum_{p|d}\log p\ll N^{-{1\over10}}
	\]
	Plugging this result back, we obtain $R_1\ll N^{9\over10}\log^2N$. Finally, combining \eqref{eqn-m10} and \eqref{eqn-r0}, we conclude that for large $N$
	\begin{align*}
		\Omega
		&\le M_1+R+N^{11\over12}\le3.94033\mathfrak S(N){N\over\log^2N} \\
		&+O\left({N\over\log^AN}+N^{9\over10}\log^2N+N^{11\over12}\right)
	\end{align*}
	This gives \eqref{eqn-omegabound}, so the proof of $\{1,2\}$ is now complete.
	\bibliographystyle{plain}
	\bibliography{refs.bib}

\end{document}